\documentclass[12pt, a4paper, reqno]{amsart}
\usepackage{amsthm, mathtools, amssymb, amsfonts, stmaryrd, latexsym, mathrsfs, todonotes, eufrak, color}
\usepackage{graphicx}
\usepackage[english]{babel}
\usepackage[all]{xy}
\usepackage{nicefrac}
\usepackage{mparhack}
\usepackage[normalem]{ulem}
\usepackage{color}

\usepackage{fancyhdr}
\usepackage{enumitem}
\usepackage{algorithmic}
\usepackage{algorithm}
\usepackage{hyperref}

\usepackage[foot]{amsaddr}

\theoremstyle{plain}
\newtheorem{theo}{Theorem}[section]
\newtheorem{lemma}[theo]{Lemma}
\newtheorem{prop}[theo]{Proposition}
\newtheorem{corollary}[theo]{Corollary}
\newtheorem{remark}[theo]{Remark}
\newtheorem{example}[]{Example}

\theoremstyle{definition}
\newtheorem{definition}[]{Definition}

\newcommand{\F}{\mathbb{F}}

\newcommand{\Z}{\mathbb{Z}}
\newcommand{\Q}{\mathbb{Q}}
\newcommand{\R}{\mathbb{R}}

\newcommand{\N}{\mathbb{N}}

\newcommand{\Ann}{{\rm Ann}}

\renewcommand{\epsilon}{\varepsilon}

\newcounter{enumi_saved}

\makeatletter
\def\imod#1{\allowbreak\mkern10mu({\operator@font mod}\,\,#1)}
\makeatother

\allowdisplaybreaks

  \title{On Fuchs' Problem about the group of units of a ring}
  
\begin{document}

\maketitle
\begin{center}
\author{\textsc{{Ilaria Del Corso}\footnote{(I. Del Corso): Dipartimento di Matematica Universit\`a di Pisa, e-mail: ilaria.delcorso@unipi.it} }\author{ \textsc{and  Roberto Dvornicich}\footnote
{(R. Dvornicich): Dipartimento di Matematica Universit\`a di Pisa, e-mail: roberto.dvornicich@unipi.it}  }}
\end{center}

\begin{abstract}
In \cite[Problem 72]{Fuchs60} Fuchs posed the problem of characterizing  the  groups which are the groups of units of  commutative rings. 
In the following years, some partial answers have been given to this question in particular cases. In a previous paper \cite{DDcharp} we dealt with finite characteristic  rings. In this paper we consider Fuchs' question for finite groups and we address this problem in two  cases. Firstly, we study the case of torson-free rings and we obtain a complete classification of the finite groups of units which arise in this case. Secondly, we examine the case of characteristic zero rings obtaining, a pretty good description of the possible groups of units equipped with families of examples of both realizable and non-realizable groups.
The main tools to deal with this general case are the Pearson and Schneider splitting of a ring \cite{PearsonSchneider70}, our previous results on finite characteristic rings \cite{DDcharp} and our classification of the groups of units of torsion-free rings.

As a consequence of our results we completely answer Ditor's  question \cite{ditor} on the possible cardinalities of the group of units of a ring. 
\end{abstract}

MSC: 16U60,13A99, 20K15.

\section{Introduction}

In this paper we consider the famous problem posed by Fuchs in \cite[Problem 72]{Fuchs60}: 

\begin{itemize}
\item[]{\sl Characterize the groups which are the groups of all units in a commutative and associative ring with identity.}
\end{itemize}

A partial approach to this problem was suggested  by Ditor \cite{ditor} in 1971, with the following less general question: 
\begin{itemize}
\item[]{\sl Which whole numbers can be the number of units of a ring? }
\end{itemize}

In the following years, these questions  have been considered by many authors.
A first result is due to  Gilmer \cite{Gilmer63}, who considered the case of {\it finite commutative rings}, classifying  the possible cyclic groups that arise in this case. 
A contribution towards the solution of the problem can be derived from the results by  Hallett and Hirsch \cite{HallettHirsch65}, and subsequently by Hirsch and Zassenhaus
\cite{HirschZassenhaus66}, combined with  \cite{Corner63}. In fact, Hallett and Hirsch dealt with the different  problem of describing the group $Aut(G)$ when $G$ is a torsion free group, and found necessary conditions for a finite group to be  the automorphism group of a torsion free group.
Now, if $G$ is an abelian group, then $Aut(G)$ is the group of units of the ring
$End(G)$ and a deep result of Corner  \cite{Corner63} shows that any {\it countable, torsion-free and reduced ring}
is in fact isomorphic to $End(G)$ for some countable and torsion-free abelian 
group $G$. It follows that the results in  \cite{HallettHirsch65}
produce  a proof (although rather indirect) that 
if  a finite group is  the group of units  of a  reduced and torsion free ring, then it must satisfy some necessary conditions, namely, it must be a subgroup of a direct product of groups of a given  family.  

Later on, Pearson and Schneider \cite{PearsonSchneider70} combined the result of Gilmer and the result 
of Hallett and Hirsch to describe explicitly all possible {\it finite cyclic groups} that can occur as $A^*$ for a commutative ring $A$. 

Other instances on this subject were considered also by other authors,
like Eldridge and Fischer \cite{EldridgeFischer67}, who considered
rings with the descending chain condition, or Dol${\rm \check{z}}$an
\cite{Dolzan02}, who classified the finite rings $A$ for which $A^*$
is a finite group whose order is either a prime power or not divisible
by 4. Recently, Chebolu and Lockridge \cite{CheboluLockridge15} \cite{CheboluLockridge17} have
studied Fuchs' problem in a different setting, and have been able to
classify the {\it indecomposable abelian groups}  and the {\it finite dihedral groups} that are realizable
as the group of units of a ring.
\smallskip

For short, we shall call {\it realizable} a finite abelian group which is the group of units of some commutative ring.

\smallskip

In a previous paper \cite{DDcharp} we considered Fuchs' and Ditor's questions in the case when $A$ is a {\it finite characteristic ring}, obtaining necessary conditions (\cite[Thm 3.1]{DDcharp}) for a group to be realizable in this case, and therefore to produce infinite families of non-realizable groups. 
On the other hand, we also  gave positive answers; in fact,  we  exhibited  large classes of groups that are realizable. Moreover, our results allowed us to completely answer Ditor's question for finite characteristic rings.

\smallskip
In this paper we consider Fuchs' question for finite groups and rings of any characteristic.
In Section \ref{domains} we classify the possible groups of units in the case of integral domains, proving the following (see Theorem \ref{domini})

\smallskip

{\bf Theorem A}
{\sl The finite abelian groups that occur as group of units of an integral domain $A$ are:

i) the multiplicative groups of the finite fields if ${\rm char}(A)>0$;
 
ii) the cyclic groups of order 2,4, or 6 if ${\rm char}(A)=0$.}

\smallskip

In Section \ref{torfree} we study the case of torson-free rings and we obtain a complete classification of the finite groups of units which arise in this case. 
The result is the following (see Theorem \ref{2a3b}) 

\smallskip

{\bf Theorem B}
{\sl 
  The finite abelian groups which are the group of units of a {torsion-free} ring $A$, 
  %\old{without nilpotents and torsion-free} 
  are all those of the form
$$(\Z/2\Z)^a\times (\Z/4\Z)^b\times (\Z/3\Z)^c$$
where $a,b,c\in\N$, $a+b\ge1$ and $a\ge1$ if $c\ge1.$ 

In particular, the possible values of $|A^*|$ are the integers $2^d3^c$ with $d\ge1$.
}
\smallskip

The more relevant part of the proof consists in producing explicit examples of rings realizing these groups of units.
In fact, for $c=0$ one can simply choose $A=\Z^{a}\times\Z[i]^b$, whereas, for $c\ge1,$ one needs a more sophisticated construction, as described in Proposition \ref{An}.
\smallskip

In Section \ref{zero} we consider  Fuchs'
question in the general case when $A$ is a characteristic zero ring and $A^*$ is finite. 
Defining $\epsilon(A)$ as the minimum exponent of 2 in the decomposition of the 2-Sylow of $A^*$ as direct sum of cyclic groups,
we prove (see Theorem \ref{2n})

{\bf Theorem C}
{\sl 
The finite abelian groups which are the group of units of a ring $A$ of characteristic 0  have the form
\begin{equation}
\label{a*}
\Z/2^\epsilon\Z\times H
\end{equation}
where $\epsilon=\epsilon(A)=1,2$ and $H$ is an abelian group.
}
\smallskip

It is known that all groups of the form $\Z/2\Z\times H$, where $H$ is any abelian group, are realizable; on the other hand, we show that, in the case when  $\epsilon(A)=2$, this is no longer true and a more subtle analysis  is necessary. 
A fundamental tool for the study of this general case is Proposition \ref{ps} (due to Pearson and Schneider) which allows to reduce the study to finite rings and to rings  whose torsion elements are nilpotent (which we call {\it type 2 } rings). Finite rings were studied in \cite{DDcharp}, whereas the study of type 2 rings is done  by using the results of Section \ref{torfree} on torsion-free rings.
The outcome is a pretty good description of the situation, which
allows us to find families of non-realizable groups (see Proposition \ref{quadrato}) and to produce large families of realizable groups (see Proposition \ref{H2}).  In particular, we can easily reobtain the classification of the cyclic realizable groups. 
 
Finally, we include an Appendix where we present some results concerning the
densities of the cardinalities of realizable groups $A^*$ in
the set of the natural numbers. In Proposition \ref{densitatotale} we prove that the density of the cardinalities of all realizable groups 
 is equal to $\frac 12$, whereas
  in Proposition \ref{density-nonilp}  we show that restricting to reduced rings  the density is zero.
 
\section{Notation and preliminary results}

Let $A$ be a ring with 1 and, as usual,
denote by $A^*$ its multiplicative group, which we shall assume to be 
{\it finite and abelian} throughout the paper. The following remark shows that for the study of all possible groups $A^*$ we may restrict to a particular class of rings.  
\begin{remark}
  \label{ZA} {\rm Let $A_0$ be the fundamental subring of $A$ (namely $\Z$ or $\Z/n\Z$ depending on whether the characteristic of $A$ is 0 or $n$) and let
    $B=A_0[A^*]$. Then trivially $B^*=A^*$.  Hence, without loss of generality, we
    shall  assume that $A=A_0[A^*]$. 
    
     It follows that for our purposes  we can
    assume that $A$ is {\em commutative, finitely generated and integral over $A_0$.} We will make this assumption throughout the paper
       }
 \end{remark}

For a commutative ring $A$, we denote by $\mathfrak{N}=\mathfrak{N}(A)$ its nilradical, namely, the ideal of the nilpotent elements. 
A ring is called reduced if its nilradical is trivial.

We note that replacing  $A$ with $A_0[A^*]$ does not change the property of being reduced because they  have the same nilradical. In fact, the nilradical of $A_0[A^*]$ is obviously contained in the nilradical of $A$; on the other hand, if $a\in A$ is such that $a^n=0$ for some $n\in\N$, then $1+a\in A^*$, so $a\in A_0[A^*]$. 
%In particular, $A$ is reduced if and only if $A_0[A^*]$ is reduced. 

An element $a\in A$ is a torsion element if its additive order is finite. $A$ is torsion-free if 0 is its only torsion element.  

For each $n\ge2$ we will denote by $\zeta_n$ a complex primitive $n$-th root of unity.
\begin{lemma}
\label{nimplicat}
{\sl
Let $A$ be  a commutative ring with 1 and assume that $A^*$ is finite. Let $x\in \mathfrak{N}$, then $x$ is a torsion element.
In particular, if $A$ is torsion-free, then it is reduced.
}
\end{lemma} 
\begin{proof}
We know that $1+ \mathfrak{N}\subseteq A^*$; this yields that  $1+ \mathfrak{N}$, and hence $\mathfrak{N}$, is finite. Let $x\in \mathfrak{N};$  the elements of the $\Z$-module $\{nx\}_{n\in\Z}$  are all nilpotents, hence they can not be all distinct. This means that $x$ is a torsion element.
\end{proof}

In the following we shall make repeated use of the following proposition (see \cite[Prop.2.2]{DDcharp}):
\begin{prop}
\label{successioneesatta}
{\sl
Let $A$ be a commutative ring and let $\mathfrak{N}$ be its nilradical. Then the sequence
\begin{equation}
\label{success}
1\to 1+\mathfrak{N}\hookrightarrow A^*\stackrel{\phi}{\to}\left(A/\mathfrak{N}\right)^*\to 1,
\end{equation}
where $\phi(x)=x+\mathfrak{N}$, is exact.
}
\end{prop}
%\begin{proof} Clearly  
%$$ 1\to 1+\mathfrak{N}\to A^*\to{A^*}/{(1+\mathfrak{N})}\to 1$$
%is exact. Moreover, the map $\Phi : A^*\to (A/\mathfrak{N})^*$  defined by $x\mapsto x+\mathfrak{N}$ is a homomorphism with kernel equal to $1+\mathfrak{N}$. Further, if $x+\mathfrak{N}\in (A/\mathfrak{N})^*$ and $(x+\mathfrak{N})( y+\mathfrak{N})=1+\mathfrak{N}$ then $xy\in1+\mathfrak{N}\subseteq A^*$, so $x\in A^*$ and $\Phi$ is surjective, so $A^*/(1+\mathfrak{N})\cong( A/\mathfrak{N})^*.$
%\end{proof}

Although  for finite characteristic rings the exact sequence  \eqref{success} always splits (see \cite[Thm 3.1]{DDcharp}), this is no longer true in general, as shown in Section \ref{zero} (see Examples \ref{nonsplit}).

\section{Integral Domains}\label{domains}

In this section we characterize the finite abelian groups $A^*$ when $A$ is a domain. 	
The fundamental result to deal with this case is 
%\new{however} 
the
famous Dirichlet's units theorem. For the sake of completeness, we
quote here the generalization of the classical Dirichlet's units
theorem to orders of a number field (see \cite[Ch.1,
\textsection12]{Neukirch99} for a proof).

\begin{theo}[Dirichlet]
\label{dirichlet}
{\sl
Let $K$ be a finite extension of $\Q$ such that $[K:\Q]=n$ and denote by $R$ an order of the algebraic integers of $K$.
Assume that among  the $n$ embeddings of $K$ in $\bar\Q$, $r$ are real (namely map $K$ into $\R$) and $2s$ are non-real ($n=r+2s$). Then
$$R^*\cong T\times\Z^{r+s-1}$$
where $T$ is the group of the roots of unity contained in $R$.
}
\end{theo}

\begin{theo}
\label{domini}
{\sl
The finite abelian groups that occur as group of units of an integral domain $A$ are:

i) the multiplicative groups of the finite fields if ${\rm char}(A)>0$;
 
ii) the cyclic groups of order 2,4, or 6 if ${\rm char}(A)=0$.
}
\end{theo}

\begin{proof}
If $A$ is a domain, then clearly also $A_0[A^*]$ is a domain, so we may assume $A=A_0[A^*]$.

If the characteristic of $A$ is a prime number $p$, then $A$ is a finite domain, hence $A$ is a field. 
It follows that $A^*$ is the multiplicative group of a finite field. Trivially, any finite field 
$\F_q$ produces its multiplicative group $\F_q^*$. 

If ${\rm char}(A)=0$, then $A$ is a finitely generated $\Z$-module, so
its quotient field $K$ is a number field and, by Remark \ref{ZA}, $A$
is contained in the ring of integers ${\mathcal O}_K$ of $K$. This guarantees that $A$ is an order of $K$ and the structure of its group of  units is described by Dirichlet's units theorem.  In particular, if $A^*$ is finite, then $r+s-1=0$, so the quotient field $K$ of $A=\Z[A^*]$ is  $\Q$ or a quadratic imaginary field. 
It follows that, if $A$ is a domain of characteristic 0 and $A^*$ is finite, then $A^*$ is 
%the possible finite abelian groups which can be obtained as of units of a domain of characteristic zero are 
a subgroup of the cyclic group of roots of unity contained in an imaginary quadratic field. 
%\old{As it is well known,  a primitive $n$-th root of unity $\zeta_n$ has degree $\phi(n)$ over $\Q$ (here $\phi$ denote the totient Euler function), and $\phi(n)\le2$ if and only if $n=1,2,3,4,6$. Since  $\pm1\in\Q$ we have the  group of units of an imaginary quadratic field can only be}
%\new{hence its group of units is a subgroup of } $\{\pm1\}$, $\langle \zeta_4\rangle$ or $\langle\zeta_6\rangle$.   \old{Finally}
Now, since $A^*$ has even  order ($\Z\subseteq A$),   the possibilities for $A^*$ are $\{\pm1\}$, $\langle\zeta_4\rangle$ and $\langle\zeta_6\rangle$. On the other hand, it is clear that  each of those groups does occur as $A^*$ for  some domain $A$; for example, we can take $A=\Z$, $\Z[\zeta_4]$ and $\Z[\zeta_6]$, respectively.
\end{proof}

 \section{{Torsion-free rings} }
 \label{torfree}
 
The case of  torsion-free rings has already been considered in the literature. 
Actually, if $A^*$ is finite and $A$ is torsion-free then, by Lemma \ref{nimplicat}, it is also reduced, and  one can deduce from \cite{HallettHirsch65} (by using the deep  result of Corner \cite{Corner63}), that, if $A^*$ is finite and abelian, then it is a
subgroup of a group of the form $(\Z/2\Z)^a\times(\Z/4\Z)^b\times(\Z/3\Z)^c$. However, the question of
which groups actually occur remained open:  
%
% in the next theorem we show that it is necessary to add a condition on the exponents  $a,b,c$ and this allows to completely
%classify the groups of units in this case. 
%
in the next theorem 
 we find  a condition on the exponents $a,b,c$ which is necessary and sufficient  for the group to be realizable, so we  completely
classify the groups of units in this case. 

We present an independent proof of the necessity of the condition. As to the sufficiency,  we produce examples of all realizable groups by means of the rather sophisticated construction given in Proposition \ref{An}.

\begin{theo}
  \label{2a3b} {\sl 
  The finite abelian groups which are the group of units of a {torsion-free} ring $A$, 
  %\old{without nilpotents and torsion-free} 
  are all those of the form
$$(\Z/2\Z)^a\times (\Z/4\Z)^b\times (\Z/3\Z)^c$$
where $a,b,c\in\N$, $a+b\ge1$ and $a\ge1$ if $c\ge1.$ 

In particular, the possible values of $|A^*|$ are the integers $2^d3^c$ with $d\ge1$.
}
\end{theo}
 In the following lemmas we assume $A$ to be a torsion-free ring such that $A^*$ is finite.
  \begin{lemma}
 \label{lemma1}
 {\sl Let  $\alpha\in A^*$ and let $\varphi_\alpha\,:\, \Z[x]\to A$ be the  homomorphism defined by $p(x)\mapsto p(\alpha)$.
Then  $\ker(\varphi_\alpha)$ is a principal ideal generated by a non-constant  polynomial.
}
 \end{lemma}
 
 \begin{proof}
Let $I=\ker(\varphi_\alpha)$ and let  $I=(f_1(x),\dots,f_n(x))$.
% for some $f_1(x),\dots,f_n(x)\in Z[x].$  The extension of $I$ to $\Q[x]$ is a principal ideal, 
Denote by $I\Q[x]$  the extension of $I$ to $\Q[x]$. Then,  $I\Q[x]=(d(x))$, where  
\begin{equation}
\label{dx}
d(x)=\sum_{i=1}^na_i(x)f_i(x)
\end{equation}
 for some  $a_i(x)\in\Q[x]$.
% ; we can also assume  
 Up to multiplication by a non-zero constant (this does not change the ideal generated by $d(x)$ in $\Q[x]$), in equation \eqref{dx} we can assume that 
 %$d(x)\in\Z[x]$ and $d(x)$ primitive.
 $a_i(x)\in\Z[x]$ for $i=1,\dots, n$, and therefore  $d(x)\in I$. Let $d(x)=c d_1(x)$, where $c\in\Z\setminus\{0\}$ and $d_1\in\Z[x]$ is primitive. Since $d(\alpha)=cd_1(\alpha)=0$ and $A$ is torsion-free, we have that
 %; with this choice of $d(x)$, clearly, $d(x)\Z[x]\subseteq I$. We want to show that the two ideals are equal.  
% From B\'ezout identity,  it follows that there exists $c\in\Z$ such that $cd(x)\in I$. But then $cd(\alpha)=0$ and, since $A$ is torsion-free, 
$d_1(\alpha)=0$, so $d_1(x)\in I$. On the other hand, if $f\in I$, then $d_1(x)$ divides $f(x)$ in $\Q[x]$ and, by Gau\ss\ Lemma,  $d_1(x)$ divides $f(x)$ in $\Z[x]$ too. In conclusion $I=(d_1(x))$.
 
 To conclude the proof, we note that $d_1(x)$ has positive degree, since $A$ is torsion-free
and therefore $I\cap\Z=\{0\}$ . 
 \end{proof}

 \begin{remark}
 {\rm
The previous lemma  is not true if the ring $A$ is not torsion-free. In fact, let $A= \Z[x]/I$ where
 $I=(x^2,2x)$ and let    $\alpha =\bar x \in A$; then $I$  is the kernel of the map  $p(x)\mapsto p(\alpha)$  but it is not principal.
 }
 \end{remark}

 \begin{lemma}
 \label{lemmapk}
 {\sl
 Let $p$ be a prime, let $k>0$ and assume that $p^k\ne 2,3,4.$ Then in $A^*$ there are no elements of order $p^k$. 
 }
 \end{lemma}
  \begin{proof}
 Assume, on the contrary, that there exists a ring $A$ fulfilling the hypotheses  and 
 an element $\alpha\in A^*$ of order $p^k$,  with  $p^k\ne 2,3,4$. Note that in this case $\phi(p^k)>2$.
 By Lemma \ref{lemma1}, we know that $\Z[x]/\ker(\varphi_\alpha)$ is a subring of $A$ and $\ker(\varphi_\alpha)$ is a principal ideal generated by a non constant polynomial $d(x)\in\Z[x]$, say. We claim that $\Z[x]/(d(x))$ contains a unit of infinite order.

%Assume, on the contrary, 
%that there exists a ring $A$ fulfilling the hypotheses of the proposition and 
% an element $\alpha\in A^*$ of order $p^k$, where $p$ is a prime, $k\ge1$ and $p^k\ne 2,3,4$. Note that in this case $\phi(p^k)>2$.
% 
% Let $\varphi\,:\, \Z[x]\to\Z[\alpha]\subseteq A$ be
%  the homomorphism induced by $x\mapsto \alpha$ and
%  let $I=\ker(\varphi)$.
%%  ; then $I\supseteq (x^p-1)$ and  
%%  $I\cap\Z=\{0\}$. 
%We show that $I$ is a principal ideal generated by a non constant polynomial.
%  % generated by a non-constant divisor $f(x)$ of $x^p-1$. 
%  In fact, $I\Q[x]=(d(x))$ where  $d(x)$ is a gcd of a set of generators of $I$: note that $d(x)$ has positve degree since $I\cap\Z=\{0\}$ and  that it can be chosen integral and primitive. From B\'ezout identity \todo{spiegare??} it follows that there exists $c\in\Z$ such that $cd(x)\in I$. But then $cd(\alpha)=0$ and, since $A$ is torsion-free, $d(\alpha)=0$ so $d(x)\in I$. On the other hand, if $f\in I$ then $d(x)$ divides $f(x)$ in $\Q[x]$ and, by Gau\ss\ Lemma,  $d(x)$ divides $f(x)$ in $\Z[x]$ too. In conclusion $I=(d(x))$.
%  
 The element $\alpha$ has order $p^k$,  so $\ker(\varphi_\alpha)=(d(x))\supseteq (x^{p^k}-1)$, namely,
 $$d(x)\mid (x^{p^k}-1)=(x^{p^{k-1}}-1)\Phi_{p^k}(x),$$
 where $\Phi_{p^k}(x)$ denotes the $p^k$-th cyclotomic polynomial.  
  Since $\Phi_{p^k}(x)$ is irreducible in $\Z[x]$ we have two cases: either $(d(x),\Phi_{p^k}(x))=1$ or $(d(x),\Phi_{p^k}(x))=\Phi_{p^k}(x)$.
  
  However, if $(d(x),\Phi_{p^k}(x))=1$, then $d(x)\mid x^{p^{k-1}}-1$,
  hence $\alpha^{p^{k-1}}=1$, contrary to our assumption.
  
  It follows that necessarily $(d(x),\Phi_{p^k}(x))=\Phi_{p^k}(x)$. Set $d(x)=h(x)\Phi_{p^k}(x)$ for some polynomial $h(x)\in \Z[x]$. 
  The Chinese Remainder Theorem gives an injection
  \begin{equation}
  \label{injection}
  \psi\,  :\, \Z[x]/(d(x)) \to \Z[x]/(h(x)) \times \Z[x]/(\Phi_{p^k}(x)). 
  \end{equation}
  Clearly the composition of the map $\psi$ with the projection to each factor is surjective, whereas
  $Im(\psi)=\{(\bar a,\bar b)\mid a-b\in(h(x),\Phi_{p^k}(x))\},$
   where $\bar a$ and $\bar
  b$ denote the projections of elements of $\Z[x]$ onto their
  respective quotients. 
  
  We note that $(h(x),\Phi_{p^k}(x)) =(p,h(x)):$ in fact, one sees immediately that 
  $\Phi_{p^k}(x)\equiv p\pmod{x^{p^{k-1}}-1},$  so $p \in (x^{p^{k-1}}-1,
  \Phi_{p^k}(x))\subseteq (h(x),\Phi_{p^k}(x))$. On the other hand,
  $\Phi_{p^k}(x)\equiv {(x^{p^{k-1}}-1)^{p-1}}\pmod{p}$, hence
  $\Phi_{p^k}(x)\in(p,x^{p^{k-1}}-1)\subseteq (p,h(x))$. 
    Summing up we have
$$Im(\psi)=\{(\bar a,\bar b)\mid a-b\in(p,h(x))\}.$$ 
  
We shall obtain a contradiction by showing that $Im(\psi)$, and hence
$ \Z[x]/(d(x))$, contains a unit of infinite order.  Since
$\phi(p^k)>2$, by Dirichlet's unit Theorem (Thm.\ref{dirichlet})
$\Z[\zeta_{p^k}]\cong \Z[x]/(\Phi_{p^k}(x))$ contains a unit
$\epsilon$ of infinite order.  Let $g(x)\in\Z[x]$ be a representative
of $\epsilon$, {\it i.e.}, $\epsilon=g(\zeta_{p^k})$; clearly, since 
$g(\zeta_{p^k})$ is a unit, it does not belong to the prime ideal
$(p,\zeta_{p^k}-1)$, hence $g(\zeta_{p^k})\equiv g(1)\not\equiv
0\pmod{(p,\zeta_{p^k}-1)}$ and therefore $g(1)\not\equiv 0 \pmod p$.

  It follows that $g(x)^{\phi(p^k)}\equiv g(x^{p^{k-1}})^{p-1}\equiv g(1)^{p-1}\equiv 1 \pmod {(p, x^{p^{k-1}}-1)}$, and, since $h(x)$ divides  $x^{p^{k-1}}-1$, the same congruence holds also modulo $(p,h(x))$. Hence  
 $(\bar 1,\overline {g(x)}^{\phi(p^k)})\in Im(\psi)$ is a unit of infinite order, a contradiction.
\smallskip

 \end{proof}
 
  \begin{lemma}
 \label{lemma12}
 {\sl
 A does not contain a unit $\alpha$ of multiplicative order 12 such that $\alpha^6=-1$.
   }
 \end{lemma}
  \begin{proof}
 Assume, on the contrary, that there exists an element $\alpha\in A^*$ such that $\alpha^6=-1$ and 
consider the substitution homomorphism  $\varphi_\alpha\colon \Z[x]\to A$ of Lemma \ref{lemma1}. 
Let ${\rm ker}(\varphi_\alpha)=(d(x))$; we have that $d(x)|x^6+1,$ hence $d(x)=x^2+1, x^4-x^2+1,$ or $x^6+1$. 
The first two cases must be excluded; in fact, the first one would imply  $\alpha^4=1$ and the second case would imply that $\Z[x]/\ker(\varphi_\alpha)\cong\Z[\zeta_{12}]$ which contains units of infinite order (see Thm. \ref{dirichlet}).
So we are left to examine the case $d(x)=x^6+1$. We consider again the injection given by the Chinese Remainder Theorem
$$
\psi\colon\Z[x]/(x^6+1) \to \Z[x]/(x^2+1) \times \Z[x]/(x^4-x^2+1).
$$
%and the canonical ring isomorphism
%$$\theta\colon\Z[x]/(x^2+1) \times \Z[x]/(x^4-x^2+1) \cong\Z[i]\times\Z[\zeta_{12}].$$
Arguing as in the proof of Lemma \ref{lemmapk}, we have 
$$Im(\psi)=\{(\overline{a(x)},\overline{b(x)})\mid a(x)-b(x)\in (x^2+1,x^4-x^2+1)=(3,x^2+1)\},$$
%  and
%$$Im(\theta\circ\psi)=\{({a(i)},b(\zeta_{12}))\mid a(i)-b(\zeta_{12})\in (3,-\zeta_3+1)\}.$$  
we show that $Im(\psi)$ contains a unit of infinite order.

 Clearly, $\Z[x]/(x^4-x^2+1) \cong\Z[\zeta_{12}]$ so it contains a unit $\overline{g(x)}$ of infinite order, which, in particular, does not belong to any proper ideal of $\Z[x]/(x^4-x^2+1)$. 
This gives that any representative $g(x)$ of the class $\overline{g(x)}$ does not belong to the proper ideal $( 3, {x^2+1})$ of  $\Z[x]$. Since $\Z[x]/(3, x^2+1)\cong \F_9$, we have that $1-g(x)^8\in (3, x^2+1)$, therefore, $(\bar1, \overline{g(x)}^8)\in Im(\psi)$ and is a unit of infinite order,  contradicting the finiteness of $A^*$.

%We show that $Im(\theta\circ\psi)$ contains a unit of infinite order.
%
%
%
%Consider the ring  $\Z[\zeta_{12}]$: by \new{TEOREMA} there exists only one prime $P$ containing 3, namely, $P=(\zeta_3-1)$.
%%It results that $P\supsetneq(3)$ and $\Z[\zeta_{12}]/P\cong\F_9$. 
%%Moreover, $P\cap\Z[i]=(3)$ and $\Z[i]/(3)\cong\F_9$.
%
%Let $g(x)\in\Z[x]$ be such that $g(\zeta_{12})$ is a unit of infinite order of $\Z[\zeta_{12}]$. 
%We have $ g(\zeta_{12})\not\in P$, hence $g(\zeta_{12})^8\equiv 1\pmod P$, since  $\Z[i]/(3)\cong\F_9$. It follow that  
%%$$g(i)\equiv g(\zeta_{12}^3)\equiv g(\zeta_{12})^3\not\equiv 0\pmod P,$$
%%so $g(i)^8\equiv 1\pmod 3$ and $g(\zeta_{12})^8\equiv 1\pmod P$.
%%Now, $1-g(\zeta_{12})^8 =(1-g(i)^8)+(g(i)^8-g(\zeta_{12})^8)\in P$ and hence  
%$(1, g(\zeta_{12})^8)\in Im(\theta\circ\psi)$ and is a unit of infinite order,  contradicting the finiteness of $A^*$.
 \end{proof}
 \begin{remark}
 {\rm
 The condition $\alpha^6=-1$ of the lemma can not be relaxed to $\alpha^{12}=1$; for example the ring $\Z[\zeta_3]\times\Z[i]$ is a torsion free ring, has a finite group of units and $\alpha=(\zeta_3, i)$ has order 12.

 }
 \end{remark} 

\smallskip

\begin{proof} [Proof of Theorem \ref{2a3b}]
Lemma \ref{lemmapk} ensures that the cyclic factors of $A^*$ can have only order 2, 3 or 4. It follows that  
$$A^*\cong(\Z/2\Z)^a\times (\Z/4\Z)^b\times (\Z/3\Z)^c$$
for some $a,b,c\ge0$. Moreover, $\Z^*\subseteq A^*$ so $2$ divides $|A^*|$ and $a+b\ge1$. 
To show that, if $c\ge1$, then $a\ge1$, we use Lemma \ref{lemma12}.
%This proves that $A^*$ is of the form $(\Z/2\Z)^a\times (\Z/4\Z)^b\times (\Z/3\Z)^c$; we now show that the conditions on the exponents are necessary. Clearly, $a+b\ge1$ since $-1\in A^*$. Let now suppose $c\ge1$ and 
Assume, by contradiction, that $a=0$;
in this case $-1$ is the square of an element of $\alpha\in A^*$, and since $c\ge1$ there exist $\gamma\in A^*$ of order 3. It follows that  $(\alpha\gamma)^6=-1$: this gives a contradiction by Lemma \ref{lemma12}.

To conclude the proof we have to show that all abelian groups of the form 
$$(\Z/2\Z)^a\times (\Z/4\Z)^b\times (\Z/3\Z)^c,$$
 with $a+b\ge1$ and $a\ge1$ if $c\ge1$, actually occur as  groups of units of suitable rings of characteristic 0, without nilpotents and torsion free. 
In fact, for $c=0$ one can choose $A=\Z^{a}\times\Z[i]^b$; for $c\ge1$ the following proposition provides a family of rings $\{A_n\}_n$ such that  
$A=\Z^{a-1}\times\Z[i]^b\times A_{c-1}$ has group of units of the required form.
\end{proof}

\begin{prop}
\label{An}
{\sl For $n\ge 0$ let 
$$
A_n = \Z[x,y_1,\dots,y_n]/(x^2+x+1,\{y_i^2+y_i+1\}_{i=1,\dots,n}).  
$$
Then $A_n$ is a ring without nilpotents and torsion free and $$A_n^*\cong\Z/2\Z\times\left(\Z/3\Z\right)^{n+1}.$$ }
\end{prop}

\begin{proof}
Let $I=(x^2+x+1,\{y_i^2+y_i+1\}_{i=1,\dots,n})$.  
For each $S\subseteq\{1,\dots,n\}$ let 
$${\mathbf y}_S=\prod_{i\in S} y_i;$$
then $\{{\mathbf y}_S+I,x{\mathbf y}_S+I\}_{S\subseteq\{1,\dots,n\}}$ is
a $\Z$-basis for $A_n$. It follows that the set $\mathcal P$ of polynomials of $\Z[x,y_1,\dots,y_n]$ of degree less than 2 in each of the variables is a complete and irredundant set of representatives of the elements of $A_n$.

For $1\le i\le n$ define $f_i=y_i-x$ and
$g_i=y_i-x^2$. Moreover, for each $S\subseteq\{1,\dots,n\}$, put
$$
I_S = I + (\{f_i\}_{i\not\in S},\{g_i\}_{i\in S}\})
$$
and 
$$B_S={\Z[x,y_1,\dots,y_n]}/{I_S}.$$
Define also
$$B=\prod_SB_S,$$
 where the product is taken over all subsets $S$ of $\{1,\dots,n\}$.
 
\begin{lemma}
\label{b}
{\sl 
For each $S\subseteq\{1,\dots, n\},$ 
$$B_S\cong\Z[\zeta_3]$$
whence,  $B\cong(\Z[\zeta_3])^{2^n}$ and
 $B^*\cong(\Z/2\Z)^{2^n}\times(\Z/3\Z)^{2^n}$.
}
\end{lemma}
\begin{proof} In $B_S$ every $y_i$ is congruent to either $x$ or to $x^2$ according to $i\not\in S$ or $i\in S$. In any case, $B_S\cong\Z[x]/(x^2+x+1) \cong \Z[\zeta_3]$ for all $S$.

The statement about $B$ is clear since $B$ is the product of $2^n$ factors all isomorphic to $\Z[\zeta_3]$. Finally, by Theorem \ref{dirichlet}, $\Z[\zeta_3]^*\cong\Z/6\Z$ and the result on $B^*$ follows. 
\end{proof}

Let $\pi_S\,:\, \Z[x,y_1,\dots,y_n]\to B_S$ be the
canonical projection  and let
$$
\pi\,:\,\Z[x,y_1,\dots,y_n]\to B.
$$
defined by $\pi(\xi)=(\pi_S(\xi))_S$ for all $\xi\in\Z[x,y_1,\dots,y_n].$

\begin{lemma}
\label{iis}
{\sl The map  $\pi$ induces an injective ring homomorphism 
$$\psi \,:\,A_n\to B.
$$ 
}
\end{lemma}  
\begin{proof}
The map $\pi$ induces an injective homomorphism 
$A_n\to B$ if and only if $\ker(\pi)=I$. By the Chinese Remainder Theorem, this is the case if and only if 
$I=\cap_SI_S $. 

Clearly, $I\subseteq\cap_SI_S $ and we are left to verify the reverse inclusion.

Let $f\in \cap_S I_S$;
 and let $r\in\mathcal P$ be its canonical representative in $A$, namely $f+I=r+I$.
% \subseteq
%\cap_S \tilde{I}_S \cap \Z[x,y_1,\dots,y_n] =\tilde{I} \cap
%\Z[x,y_1,\dots,y_n]$. Then \todo{a che serve questa formula?}
%
%$$f=h_0\cdot(x^2+x+1)+ \sum\limits_{i=1}^n
%h_i\cdot(y_i^2+y_i+1)
%$$
%with $h_0,h_i\in\bar\Q[x,y_1,\dots,y_n]$; on the other hand, 
%diving $f$
%by the monic polynomial $x^2+x+1$ we obtain $f=q_0(x^2+x+1)+r_0$ where
%$q_0,r_0\in \Z[x,y_1,\dots,y_n]$ and the degree in $x$ of $r_0$ is is strictly less than 2; similarly,
%$r_0=q_1(y_1^2+y_1+1)+r_1$ and inductively we get 
%$$
%f=q_0(x^2+x+1) +
%\sum\limits_{i=1}^n q_i(y_i^2+y_i+1) + r_n,
%$$
%where $q_i,r_n\in
%\Z[x,y_1,\dots,y_n]$ and the degree of $r_n$ in
%  each of the variables is less than 2.
  We claim that $r=0$, so $f\in I$. 
  
Denote by $\bar \Q$ an algebraic closure of $\Q$ and let  $\tilde{I}$ and $\tilde{I}_S$ be the extensions of  $I$ and
$I_S$ to the ring $\bar \Q[x,y_1,\dots,y_n]$.  
Considering the sets of zeroes
in $\bar\Q^{n+1}$ we have $V(\tilde{I})=\cup_S V(\tilde{I}_S)$ and, since the ideals
$\tilde{I}$ and $\tilde{I}_S$ are radical, this yields $\tilde{I}=\cap_{S} \tilde{I}_S$.  

Now, $ \cap_S I_S\subseteq
\cap_S \tilde{I}_S \cap \Z[x,y_1,\dots,y_n] =\tilde{I} \cap
\Z[x,y_1,\dots,y_n]$, so $f\in\tilde I$ and hence $r\in \tilde I.$ 
Considering that  $\{{\mathbf y}_S+\tilde I,x{\mathbf y}_S+\tilde I\}_{S\subseteq\{1,\dots,n\}}$ is a $\bar \Q$-basis for $\bar \Q[x,y_1,\dots,y_n]/\tilde I$, and
$$\tilde I=r+\tilde I=\sum_S (h_S{\mathbf y}_S+g_Sx{\mathbf y}_S+\tilde I),$$
we get $r=0.$

It follows that
  $\cap_S I_S \subseteq I$ and the equality holds.
\end{proof}
From Lemmas \ref{b} and \ref{iis} it follows that $A_n$ is isomorphic to a subring of $(\Z[\zeta_3])^{2^{n}}$ (in particular this shows that $A_n$ is torsion-free) and that $A_n^*$
is isomorphic to a subgroup of $(\Z/2\Z)^{2^n}\times(\Z/3\Z)^{2^n}$. 
To prove the Proposition we are left to show that the cardinality of the 2-Sylow $G_{n,2}$ and of the 3-Sylow $G_{n,3}$ of $A_n^*$ are $2$ and $3^{n+1}$, respectively.

Arguing as in the proof of Lemma \ref{lemmapk}
we can easily prove the following 
\begin{enumerate}
\item[(i)] if  $(a_S+I_S)_S\in Im (\psi)$, then $a_S-a_T\in I_S+I_T$ for all  $S\ne
T$;
\item[(ii)]  $I_S+I_T=(3,x-1,y_1-1,\dots,y_n-1)$ for all $S\ne T$.
\end{enumerate}

We show that the group $G_{n,2}$ has 2 elements. In fact, $\psi(G_{n,2})$ is a subgroup of $\{(\pm 1+I_S)_S\}$; on the other hand, if $(a_S+I_S)_S\in \psi(G_{n,2})$, by (i), $a_S-a_T\in I_S+I_T$ for all $S,T$: since $a_S=\pm1$ and $2\not \in I_S+I_T$ for all $S\ne T$, then an element of $\psi(G_{n,2})$ must have all coordinates represented by  $1$ or all represented by  $-1$, hence $G_{n,2}=\{\pm1+I\}$.
 
\smallskip

We shall now prove that $|G_{n,3}|=3^{n+1}$. 
Firstly, $G_{n,3}$ contains the set  ${\mathcal U}=\{x^{\epsilon_0}y_1^{\epsilon_1}\cdots y_n^{\epsilon_n}+I\, |\, \epsilon_i\in\{0,1,2\}\}$  and so $|G_{n,3}|\ge3^{n+1}$, since the elements of $\mathcal U$ are all distinct. To prove the last assertion it is enough to note that the canonical representatives in $\mathcal P$ of the elements of $\mathcal U$ are all distinct, since they have different factorization in $\Z[x, y_1,\dots,y_n]$.

We have to show the converse inequality. Clearly, 
$\psi(G_{n,3})$ is  a subgroup of the 3-torsion subgroup of $B^*$, namely, $\{(x^{\epsilon_S}+I_S)_S\mid
\epsilon_S=0,1,2\}$.
%By (i) and (ii), $(x^{\epsilon_S}+I_S)_S\in
%Im(\psi)$ if and only if $x^{\epsilon_S}-x^{\epsilon_T}\in (3,x-1,y_1-1,\dots,y_n-1)$, for all $S\ne
%T$. We are left with the task of characterizing the $\{\epsilon_S\}_S$ satisfying the preceding relations. 

Consider the projection 
$$pr\,:\, B\to B/3B=\prod_S B_S/3B_S,$$  let $(\bar b_S)_S=pr((b_S+I_S)_S)$  and $\bar\psi=pr\circ\psi$.
%  $\bar\psi=(\bar\psi_S)_S$. 
When restricting $pr$ to the 3-torsion elements of $B^*$ we obtain an injective map (for each $S$ the elements $1, x, x^2$ are pairwise not congruent modulo $3B_S$),
whence $|G_{n,3}|=|\bar\psi(G_{n,3})|$.

Since $B/3B\cong (\Z[\zeta_3]/3\Z[\zeta_3])^{2^n}$ and $\dim_{\F_3}\Z[\zeta_3]/3\Z[\zeta_3]=2$, then $B/3B=\prod_S B_S/3B_S$ is a vector space over $\F_3$ of dimension $2^{n+1}$. 
For each $S$ consider the base of $B_S/3B_S$ given by the classes of 1 and of $x-1$. We shall use the notation $\bar b=(\bar b_S)_S\in B/3B$ and $\bar b_S=b_{0,S}+b_{1,S}\overline{(x-1)}$ where $b_{i,S}\in \F_3$.
Note that the projection of the 3-torsion elements of $B^*$ is a subset of $B/3B$ contained in the affine subspace 
$$W=\{(\bar b_S)_S \mid  b_{0,S}= \bar1\} .$$
Put $V=\bar\psi(A_n)$; then $V$ is a vector subspace of $B/3B$, so 
$\bar\psi(G_{n,3})\subseteq V\cap W$. We shall derive an upper bound for 
$|\bar\psi(G_{n,3})|$  by bounding the dimension of the affine subspace $V\cap W$.

We now determine a set of linearly independent relations fulfilled by the elements of $V$.

 By (i) and (ii), if $v=(v_{0,S}+v_{1,S}\overline{(x-1)})_S\in V$ then 
 \begin{equation}
 \label{rel1}
 v_{0,S}=v_{0,T}\qquad \forall S\neq T.
 \end{equation}
 Considering $S=\varnothing$ and $T\neq \varnothing$ , we get $2^n-1$ independent linear relations among
the elements of $V$, whereas the other relations depend on these. 
This gives $\dim V\le 2^{n+1}-2^n+1=2^n+1$. 

Further, we show that the elements of $V$ verify the following set of relations:
\begin{equation}
\label{rel}
\sum_{T\subseteq U} (-1)^{|T|} v_T=0 
 \qquad \forall v=(v_T)_T\in V 
\end{equation}
for all $U\subseteq  \{1,\dots,n\}$  with $|U|\ge 2$.

% Let $U$ be any subset of $\{1,\dots,n\}$ with $|U|\ge 2$.
Note that $\{\bar\psi({\mathbf y}_S), \bar\psi({x\mathbf y}_S)\}_S$ is
a set of generators of $V$, hence it is enough to show that relations
\eqref{rel} are verified by these elements. Clearly
$\bar\psi_T({\mathbf y}_S)=\bar x^{|S|+|S\cap T|}$. If $U\not\subseteq S$,
then, for each $r\ge0$,  among the subsets $T$ of $U$  with $|S\cap T|=r$  one half has even cardinality and one half has  odd cardinality, so the sum clearly
vanishes. If $U\subseteq S$, and $|U|=m$ then the relation becomes
$$
\sum_{T\subseteq U} (-1)^{|T|}\bar x^{|S|+|T|}= \sum_{i=0}^m (-1)^i \binom{m}{i}\bar x^{|S|+i} = 
\bar x^{|S|}(\bar1-\bar x)^m=0, 
$$
since  $(1-x)^m\in 3B$ for $m\ge 2$. 
An analogous computation proves that \eqref{rel} holds also for $v=\bar\psi(x{\mathbf y}_S)$. 

Each of the equations \eqref{rel} splits into the two equations 
\begin{equation}
\label{rel_vi}
\sum_{T\subseteq U} (-1)^{|T|} v_{i,T}=0 
 \qquad \forall v=(v_T)_T\in V,  \quad i=0,1.
\end{equation}

The relations with $i=0$ follow directly from equations \eqref{rel1}.

The relations with $i=1$ are clearly independent from the previous ones and also among themselves (for example because one can order the subsets $U$ so that the matrix of the relations becomes triangular).
There is one relation for each subset 
$U$ of $\{1,\dots,n\}$ of cardinality at least 2, and hence we get $2^n-n-1$ relations.
 This gives 
$$
\dim V \le 2^n +1 - 2^n +n +1 = n+2.
$$
Since $V\not\subseteq W$,  the dimension of $V\cap W$ is at most $n+1$ and therefore $|\bar\psi(G_{n,3})|\leq3^{n+1}$, as desired.
\end{proof}

\section{General characteristic zero rings}
\label{zero}
In this section we consider the general case of a characteristic
  zero ring $A$. In Theorem \ref{2n} we prove that all groups of the
  form $\Z/2\Z\times H$, where $H$ is any abelian group, are
  realizable as groups of units. 
 Moreover, when  $A^*$ does not have a cyclic factor of order 2 in its decomposition,  then it must have a cyclic direct factor of order 4, and a more subtle analysis  is required. 
% However, it may happen that the
%  group $A^*$ does not have a cyclic direct factor of order 2.  If
%  this is the case, a more subtle analysis is required. 
  For this case
  we give some necessary conditions and some sufficient conditions for
  the group to be realizable (see Proposition \ref{H2}, Proposition
  \ref{quadrato} and the summary in Section \ref{sec:Concluding
    remarks}).  
%    As mentioned in the introduction, a complete
%  classification seems quite awkward, but in any case we severely
%  restrict the space between realizable and non-realizable groups. 

% Even if this is not a complete classification, for the general case. In the case when  $A^*$ does not have a cyclic factor of order 2 in its decomposition,  a more subtle analysis  allows us to find some more sufficient conditions and some more necessary conditions 

\begin{definition}
{\rm
For a ring $A$ we define $\epsilon(A)$ as the minimum exponent of 2 in the decomposition of the 2-Sylow of $A^*$ as direct sum of cyclic groups.}
\end{definition}

The following theorem describes almost precisely the group of units in the case of characteristic zero rings. The second statement of the theorem can be found also in \cite[Thm 129.1]{fuchs1970infinite};
we include the proof using our notation,
%; we rewrite the construction of the example in our notation for the sake of completeness and 
since we will need a similar construction in the proof of Proposition \ref{H2}.
However, the main result of next theorem is the fact that, for every ring $A$ of characteristic 0, necessarily $\epsilon(A)\le 2.$ 
\begin{theo}
\label{2n}
{\sl The finite abelian groups which are the group of units of a ring $A$ of characteristic 0  have the form
\begin{equation}
\label{a*}
\Z/2^\epsilon\Z\times H
\end{equation}
where $\epsilon=\epsilon(A)=1,2$ and $H$ is an abelian group.

Conversely 
%\old{For $\epsilon=1$ and } 
for each finite abelian group $H$ there exists a ring $A$ 
%\old{of type 2} 
such that 
$$
A^*\cong\Z/2\Z\times H.
$$
%\old{In particular, the possible values of $|A^*|$ are all the  even positive integers. }
}

\end{theo}
\begin{lemma}
\label{mudix}
{\sl
Let $\alpha\in A^*$ and consider the homomorphism $\varphi_\alpha\colon\Z[ x]\to A$ defined by $p(x)\mapsto p(\alpha)$.
If $\ker(\varphi_\alpha)$ contains an irreducible polynomial $\mu(x)$, then $\ker(\varphi_\alpha)=(\mu(x))$ and $\Z[\alpha]$ is a domain.
% and the composite map 
%$$\psi\colon\Z[\alpha] \hookrightarrow A\to A/{\mathfrak N}$$
%is an injection.
%map 
%$$\psi\colon\Z[\alpha]\to A/{\mathfrak N},$$
%defined by $p(\alpha)\mapsto p(\alpha)+{\mathfrak N}$,
%is an injective homomorphism of commutative rings with 1.
}
\end{lemma}
\begin{proof}
If $\mu(x)\in\ker(\varphi_\alpha)$, then $\ker(\varphi_\alpha)=(\mu(x),f_1(x),\dots, f_k(x))$ for some $f_i\in\Z[x]$.
Now, by Gau\ss\ Lemma, $\mu(x)$ is irreducible in $\Q[x]$; moreover, $1\not\in \ker(\varphi_\alpha)\Q[x]$ (in fact,  $\ker(\varphi_\alpha)\cap\Z=\{0\}$ since ${\rm char}(A)=0$),
therefore $ \ker(\varphi_\alpha)\Q[x]=(\mu(x))$ and $\mu(x)$ divides $f_i(x)$ in $\Q[x]$, for each $i$. Finally, since $\mu(x)$ is primitive, we can conclude that $\mu(x)$ divides $f_i(x)$ also in $\Z[x]$, hence $\ker(\varphi_\alpha)=(\mu(x))$ and $ \Z[\alpha]\cong\Z[x]/(\mu(x)).$

%The map $\psi\colon\Z[\alpha] \hookrightarrow A\to A/{\mathfrak N}$ obtained by composing the inclusion with the projection is an injection;  in fact, $\ker(\psi)={\mathfrak N}\cap \Z[\alpha]=\{0\}$, since $\Z[\alpha]$ is a domain.
%
%map 
%$$\psi\colon\Z[\alpha]\to A/{\mathfrak N},$$
%defined by $p(\alpha)\mapsto p(\alpha)+{\mathfrak N}$,
%is an injective homomorphism of commutative rings with 1.
%
%
%The map $\psi$  is obtained by composing the inclusion $\Z[\alpha] \hookrightarrow A$ with the projection $A\to A/{\mathfrak N}$, hence it is a homomorphism of commutative rings with 1. Moreover, $\ker(\psi)={\mathfrak N}\cap \Z[\alpha]$ is trivial since $\Z[\alpha]$ is a domain, and this proves the lemma.
% It follows that
% $A$ contains a ring isomorphic to $\Z[x]/(\mu(x))\cong \Z[\alpha].$ Moreover, the map $\Z[\alpha] \hookrightarrow A\to A/{\mathfrak N}$ obtained by composing the inclusion with the projection is injective; in fact, ${\mathfrak N}\cap \Z[\alpha]=\{0\}$ since $\Z[\alpha]$ is a domain.
\end{proof}
\begin{remark}
\label{uno}
{\rm
We note that the homomorphism $\varphi_\alpha$ induces an injection of $\Z[\alpha]$ into $A$ preserving the identity. This property will be crucial in the following since it implies that $(\Z[\alpha])^*$ is a subgroup of $A^*$.
}
\end{remark}

\begin{proof}[Proof of Theorem \ref{2n}]
Since ${\rm char}(A)=0$, then $\Z\subseteq A$ and $2=|\Z^*|$ divides
  $|A^*|$, so $\epsilon(A)\ge1$; we have to show that $\epsilon\le2.$ Assume the contrary; then $-1=\alpha^4$ for some $\alpha\in A^*$. 
Let again $\varphi_\alpha\colon\Z[ x]\to A$ be defined by $p(x)\mapsto p(\alpha)$. Then the irreducible polynomial   $\mu(x)=x^4+1$ belongs to $\ker(\varphi_\alpha)$. By Lemma \ref{mudix}, we get $\ker(\varphi_\alpha)=(\mu(x))$. It follows that 
 $A$ contains a ring isomorphic to $\Z[x]/(x^4+1)\cong \Z[\zeta_8]$ and, by Remark \ref{uno}, $A^*$  contains $\Z[\zeta_8]^*$ as a subgroup: this is not possible, since $A^*$ is finite and $\Z[\zeta_8]^*$ is infinite by Dirichlet Theorem.

As to the converse, let $H\cong\Z/a_1\Z\times\cdots\Z/a_n\Z$ where $a_1,\dots, a_n$ are positive integers. Define
$$A= \frac{\Z[x_1,\dots,x_n]\qquad}{(a_ix_i, x_ix_j )_{1\le i,j\le n}}.$$
%\old{clearly $A$ is a type 2 ring.}
Let $\mathfrak{N}$ be the nilradical of $A$; then $\mathfrak{N}=(\bar x_1,\dots,\bar x_n)$ where, as usual, $\bar x_j$ denotes the class of $x_j$ in $A$. Clearly, $A=\Z[\mathfrak{N}]$.
The exact sequence given in \eqref{success} in this case specifies to
$$ 1\to 1+\mathfrak{N}\to A^*\to\{\pm1\}\to 1,$$
 so the sequence splits and $A^*\cong \Z/2\Z\times (1+\mathfrak{N})$. 
Moreover, in this case the map $n\to 1+n$ is an isomorphism from the additive group $\mathfrak{N}$ to the multiplicative group $(1+\mathfrak{N})$, hence 
$$A^*\cong \Z/2\Z\times\mathfrak{N}\cong \Z/2\Z\times\Z/a_1\Z\times\cdots\Z/a_n\Z.$$
\end{proof}
From the previous proposition we immediately get the following
%the  answer to Ditor question \cite{ditor} in the case of characteristic 0 rings.  
\begin{corollary}
\label{corditor}
{\sl
The possible values of $|A^*|$, when $A$ is a characteristic 0 ring with finite group of units, are all the even positive integers. 
}
\end{corollary}
The last corollary and \cite[Cor. 3.4]{DDcharp}  allows to  completely answer  Ditor's question for rings of any characteristic.

\begin{corollary}
\label{ditorgen}
{\sl
The possible values of $|A^*|$, when $A$ is a ring with finite group of units, are all the even positive integers and the finite products of integers of the form $2^\lambda-1$ with $\lambda\ge 1.$ 
}
\end{corollary}
 As a particular case we reobtain the following result of Ditor \cite{ditor}
 
\begin{corollary}
\label{corcor}
{\sl
If $|A^*|=p$ is  prime, then $p=2$ or $p$ is a Mersenne prime.
}
\end{corollary}

\subsection{The case $\mathbf{\epsilon(A)=2}$.}\label{sec:epsilon2}

In Theorem \ref{2n} we classify the groups of units of the rings $A$
with $\epsilon(A)=1$, showing that any abelian group $H$ can appear in
equation \eqref{a*} in this case.  This is no longer true for rings with
$\epsilon=2$: for example, in this case we can not have $H\cong
  \Z/11\Z$, since the cyclic group $\Z/44\Z$ is not realizable (see
  \cite {PearsonSchneider70}). Many other examples can be derived from Proposition \ref{quadrato} and Corollary \ref{nocyclic}.

%  \old{ In this case we do not give a
%  complete classification of the groups of units of rings; however, }
In the following remark we {point out} an important necessary 
condition on  $A^*$ for rings with $\epsilon=2$ and in the next proposition we exhibit  a large class of groups  $H$ which occur in this case.
\begin{remark}
{\rm
Let $\epsilon(A)=2$, then the ring $A$ must contain $\Z[i]$ as a subring. In fact, in this case $-1=\alpha^2$ for some $\alpha\in A$   and so, by Lemma \ref{mudix}, we have  $\Z[i]\cong\Z[x]/(x^2+1)\hookrightarrow A.$ 

We will see that this property is crucial for showing that  not all abelian groups $H$ can occur in \eqref{a*}  in this case (see Proposition \ref{quadrato}). 
}
\end{remark}

For any finite group $H$ and for any prime number $p$ we denote by $H_p$ the $p$-Sylow subgroup of $H$. 
\begin{prop}
\label{H2}
{\sl
Let $H$ be a finite abelian group with the following properties:
\begin{enumerate}
\item[{\rm a)}] $\forall p\equiv3\ ({\rm mod}\ 4)$  the $p$-Sylow  $H_p$ is the square of a group;
\item[{\rm b)}] $H_2$ is isomorphic to {\it P} or to $\Z/4\Z\times P$
 where 
 $$P\cong \Z/2^{e_1}\Z\times\dots\times\Z/2^{e_{2r}}\Z $$ 
 with ${2\le}e_1\le\dots\le e_{2r}$ and $e_{2j}-e_{2j-1}\le1\ \forall j=1,\dots,r$.
\end{enumerate}

Then $\Z/4\Z\times H$ is the group of units of a characteristic 0 ring. 
%\old{a type 2 ring $A$}.  
    }
\end{prop}

\begin{proof}
%We show that $\Z/4\Z\times H$ is the group of units of a type 2 ring $A$  whenever $H$ has the following property:
%\begin{equation}\label{*}
%\hbox{$\forall p\equiv3\ ({\rm mod}\ 4)$ the $p$-Sylow of $H$ is of the square of a group.}
%\end{equation}
%In fact, 
To prove this proposition we need to refine the proof of Theorem \ref{2n}. 

Let $a_1,\dots, a_n$ be any finite sequence of non-zero elements of the ring $\Z[i]$.
 Similarly to above, define
$$A= \frac{\Z[i][x_1,\dots,x_n]\qquad}{(a_ix_i, x_ix_j )_{1\le i,j\le n}}.$$
%\old{clearly, $A$ is of type 2.}
By the same argument we get
$$A^*\cong \Z[i]^*\times(1+\mathfrak{N})\cong \Z/4\Z\times\mathfrak{N}\cong \Z/4\Z\times\Z[i]/(a_1)\times\cdots\Z[i]/(a_n).$$
To understand which groups can be obtained when the $a_j$'s vary, it
is enough to describe the quotients $\Z[i]/(\pi)^h$ when $\pi$ is a
prime of $\Z[i]$.  Let $(p)=(\pi)\cap\Z$ and let $N:\Z[i]\to\Z$ denote
the usual norm.  It is well known that $N(\pi)=p^2$ or $p$ according
to $p\equiv 3\pmod4$ or not; from this, and from $(2)=(1+i)^2$, we
easily get the following group isomorphisms:
\begin{equation}
\frac{\Z[i]}{(\pi)^h}\cong
\begin{cases}\Z/p^h\Z  & {\rm if\ } p\equiv1\pmod4\\
\Z/p^h\Z\times \Z/p^h\Z & {\rm if\ } p\equiv3\pmod4\\
\Z/2^k\Z\times \Z/2^k\Z & {\rm if\ }\  p=2\ {\rm and}\ h=2k\\ 
\Z/2^{k+1}\Z\times \Z/2^{k}\Z& {\rm if\ }\ p=2\ {\rm and}\ h=2k+1.
\end{cases}
\end{equation}
Now, it is clear that, by suitable choices of the $a_j$'s, one can obtain any $p$-group if $p\equiv 1\pmod4$, and
any square of a $p$-group if $p\equiv 3\pmod4$. Moreover,  the groups $P$ described in (b) are precisely the groups that can be obtained as a finite product of groups $\Z[i]/(1+i)^h$ with $h\ge4$. 
Finally, if $A$ is a ring such that  $(A^*)_2$ is isomorphic to $P$, then, putting $B=A\times \Z[i]$, we have $(B^*)_2\cong P\times\Z/4\Z.$
 
To conclude, we note that in our construction $H$ has no direct summand isomorphic to $\Z/2\Z$, hence the rings we constructed have  $\epsilon=2$.

\end{proof}

%\old{Further information on the possible groups $A^*$  when $\epsilon(A)=2$ can be obtained from an accurate study of the exact sequence 
%$$ 1\to 1+\mathfrak{N}\to A^*\to(A/\mathfrak{N})^*\to 1.$$  
%The following proposition, proved by Pearson and Schneider, provide a split }

To further investigate the structure of $A^*$ in the case  when $\epsilon(A)=2$, it is convenient to use the splitting of a ring 
proved by Pearson and Schneider, which we recall below.
%The following proposition 
%allows to obtain the information on the group of units  by combining the study of finite rings with that of rings whose torsion ideal is contained in the nilradical.
We recall that, by Remark \ref{ZA}, we can assume that our ring fulfils  the hypothesis of the proposition.

\begin{prop}{\rm (\cite[Prop.~1]{PearsonSchneider70})}
\label{ps} 
{\sl Let $A$ be a commutative ring which is finitely generated and integral over its fundamental subring. Then $A=A_1\oplus A_2$, where $A_1$ is a finite ring and the torsion ideal of $A_2$ is contained in its nilradical.}
\end{prop}

In the paper \cite{DDcharp} we studied the groups of units of finite rings, so in the following we shall concentrate on the study of the second factor of the decomposition. 

\begin{definition}
{\rm
A commutative ring is called {\it of type 2} if its torsion ideal is contained in the nilradical.}
\end{definition}

\begin{remark}\label{5.7}
  {\rm (i) A type 2 ring has characteristic 0, since otherwise 1 would
    be a torsion element and hence nilpotent.\\
 (ii) Let $A$ be a
    type 2 ring such that $A^*$ is finite. Then, by Lemma
    \ref{nimplicat}, the torsion elements of $A$ are precisely its nilpotent elements. 
    In particular, $A$ is reduced if and only if it is
    torsion free.  }
\end{remark}

% In view of the exact sequence \eqref{success}, a foundamental step is the study of rings of type 2 that are reduced rings\new{: in our case, of torsion free rings}.
%

\subsection{Type 2 rings with $\mathbf{\epsilon(A)=2}$.}
For any ring of type 2 consider the exact sequence \eqref{success} of finite groups
 $$ 1\to 1+\mathfrak{N}\to A^*\to(A/\mathfrak{N})^*\to 1$$  
and for each prime $p$  the exact sequence induced on the $p$-Sylow 
 $$ 1\to (1+\mathfrak{N})_p\to (A^*)_p\to(A/\mathfrak{N})_p^*\to 1$$
 
Since $A$ is of type 2,  the ring $B=A/\mathfrak{N}$ is torsion-free and $B^*$ is described by Theorem \ref{2a3b}, so its $p$-Sylow
$(B^*)_p$ is trivial for $p>3$.
%Moreover, also $(B^*)_3$ is trivial in this case. In fact, assume the contrary and let $\beta\in B^*$ an element of order 3.
%By Lemma \ref{mudix}, $\psi \colon\Z[i]\to B$ is an injective homomorphism of rings such that $\psi(1)=1_B$. Let $i_B=\psi(i)$; then $i_B^4=1$ and $i_B^2=-1_B$ and it is immediate to check that $i_B\beta$ is a unit of order 12 such that $(i_B\beta)^6=-1_B$, contradicting Lemma \ref{lemma12}.
%
The following proposition shows that also $(B^*)_3$ is trivial in this case.
\begin{prop}
\label{B*3}
{\sl
Let $A$ be a type 2 ring with $\mathbf{\epsilon(A)=2}$, and let $B=A/\mathfrak{N}$. Then the 3-Sylow of $B^*$ is trivial.
}
\end{prop}
\begin{proof}
Assume the contrary and let $\beta\in B^*$ an element of order 3.
By Lemma \ref{mudix}, the ring $A$ contains the domain $\Z[\alpha]$ as a subring and they have the same identity. 
Let
$\psi\colon\Z[\alpha]\to B$ be the homomorphism
 obtained by composing the inclusion $\Z[\alpha] \hookrightarrow A$ with the projection $A\to A/{\mathfrak N}$, namely, $\psi(p(\alpha))= p(\alpha)+{\mathfrak N}$. Clearly,
$\ker(\psi)={\mathfrak N}\cap \Z[\alpha]$ is trivial since $\Z[\alpha]$ is a domain and
 $\psi(1)=1+{\mathfrak N}=1_B$.
% $\psi \colon\Z[i]\to B$ is an injective homomorphism of rings such that $\psi(1)=1_B$. 
 Let $i_B=\psi(i)$; then $i_B$ has order $4$, $i_B^2=-1_B$ and it is immediate to check that $i_B\beta$ is a unit of order 12 such that $(i_B\beta)^6=-1_B$, contradicting Lemma \ref{lemma12}.
\end{proof}

Now, note that $(1+\mathfrak{N})_p=1+\mathfrak{N}_p$; in fact, these two groups have the same cardinality, so it is enough to prove that $1+\mathfrak{N}_p\subseteq  (1+\mathfrak{N})_p$:
 this can be checked by noting that,  
if $x\in \mathfrak{N}_p$ and $ k,l\ge0$ are such that $p^kx=0$ and $x^l=0$, then for $h\ge k+l$ it results $(1+x)^{p^h}=\sum_{j=0}^{p^h}\binom{p^h}{j}x^j=1$, and hence $1+x\in(1+\mathfrak{N})_p$.
% \old{Let $x\in \mathfrak{N}_p$  and let $k,l\ge0$ be such that $p^kx=0$ and $x^l=0$, then if $h$ is large (for example $h\ge k+l$) we can check easily that $(1+x)^{p^h}=\sum_{j=0}^{p^h}\binom{p^h}{j}x^j=1$, hence $1+x\in(1+\mathfrak{N})_p$.}
%Putting $B=A/\mathfrak{N}$, 
It follows that the exact sequence on the $p$-Sylow reads as
\begin{equation}
\label{sep} 
1\to 1+\mathfrak{N}_p\to (A^*)_p\to (B^*)_p\to 1.
\end{equation}
%Since $A$ is of type 2, the ring  $B$ is torsion free and reduced; by Theorem \ref{2a3b} we get that  $(B^*)_p$ is trivial for $p>3$.
and, in particular, we get
\begin{prop}
\label{nuova}
{\sl
\begin{equation}
\label{a*p}
(A^*)_p=1+\mathfrak{N}_p \ \text{for}\ p\ge3. 
\end{equation}
}
\end{prop}

%\begin{example}
%\label{nonsplit}
%{\rm
%Let $A=\Z[x,y]/(x^p-y-1, py,y^2)$ where $p=2$ or $p=3$. It is easy to see that $A$ is a type 2 ring,  $\mathfrak{N}=\mathfrak{N}_p=(\bar y)$ and $B=\Z[x]/(x^p-1)$. It follows that $1+\mathfrak{N}=1+\mathfrak{N}_p$ has exponent $p$ and $B^*\cong\Z/2\Z\times\Z/p\Z$ . It is also immediate to verify that $\bar x\in A$ is a unit of order $p^2$,  hence the exact sequence \eqref{sep} does not split.  
%}
%\end{example}

We note that the rings constructed in the proof of Proposition \ref{H2}, are in fact type 2 rings, so the result proved there holds also if we restrict to type 2 rings. In particular, 
$A^*_p$ can be any abelian $p$-group when $p\equiv1\pmod4$, and we are left to analyze the $p$-Sylow of $A^*$ for $p=2$ and $p\equiv3\pmod4$. In this last case, the following proposition gives a constraint on the structure and on the cardinality of {$1+\mathfrak{N}_p$}.

\begin{prop}
\label{quadrato}
{\sl
%\old{If $\epsilon(A)=2$ and $p\equiv3\pmod4$, then $|1+\mathfrak{N}_p|$ is a square.}
Let $A$ be a ring with $\epsilon(A)=2$ and let  $p\equiv3\pmod4$.  
For each $j\ge1$,  the quotient $(1+\mathfrak{N}_p^j)/(1+\mathfrak{N}_p^{j+1})$ has a filtration such that all its quotients   are $\F_{p^2}$-vector spaces. 

In particular  $|1+\mathfrak{N}_p|$ is a square.
}
\end{prop}
\begin{proof}
 For all $j\ge1$ the map $x\mapsto 1+x $ induces an isomorphism between $X_j=\mathfrak{N}_p^j/\mathfrak{N}_p^{j+1}$ and $(1+\mathfrak{N}_p^j)/(1+\mathfrak{N}_p^{j+1})$.

Now, $X_j$ is a $\Z[i]$-module and an abelian finite $p$-group, so for a suitable integer $r_j$ we have the filtration 
$$
X_j\supset pX_j\supset\dots\supset p^{r_j}X_j=\{0\},
$$
whose quotients $p^lX_j/p^{l+1}X_j$ are $\Z[i]/(p)\cong\F_{p^2}$-vector spaces for all $j,l.$ It follows that there exist integers $k_j$ such that $|X_j|=p^{2k_j}$, so $|1+\mathfrak{N}_p|=p^{2k}$ where $k=\sum k_j.$ 
\end{proof}

\begin{corollary}
\label{nocyclic}
{\sl
Let $\epsilon(A)=2$ and $p\equiv3\pmod4$, then $1+\mathfrak{N}_p$ is not a cyclic group.
}
\end{corollary}
\begin{proof}
All quotients of a cyclic group are cyclic, hence they can not admit a filtration whose quotients are   $\F_{p^2}$-vector spaces. 
\end{proof}

The previous proposition can not be  generalized to a result on the structure of the group; in fact, the following example shows that for $p\equiv3\pmod4$ the group $1+\mathfrak{N}_p$ is not necessarily the square of a group.

\begin{example}
{\rm
Let $p\equiv3\pmod4$ and $A=\Z[i][x]/(p^2x,px^2,x^p+px).$ 
Clearly, $\mathfrak{N}=\mathfrak{N}_p=(\bar x)$, where $\bar x$ denotes the class of $x$ in $A$, and
 as a $\Z[i]$-module
$$\mathfrak{N}_p\cong\oplus_{j=1}^{p-1}\langle \bar x^j\rangle\cong \Z[i]/(p^2)\times\left(\Z[i]/(p)\right)^{p-2}.$$
Therefore, $|1+\mathfrak{N}_p|=|\mathfrak{N}_p|=p^{2p}.$

Let $n\in\mathfrak{N}_p$, then $n$ can be written  as $n=\sum_{j=1}^{p-1}\lambda_j\bar x^j\in \mathfrak{N}_p$,
 where $\lambda_j\in\Z[i]$ and are uniquely determined modulo $p^2$ if $j=1$ and modulo $p$ if $2\le j\le p-1$.
Then $(1+n)^p=1+px(\lambda_1-\lambda_1^p)=1$ if and only if the reduction of $\lambda_1$ modulo $p$ belongs to $\F_p$; this means that there are $p^3$ possibilities for $\lambda_1$ and $p^2$ possibilities for $\lambda_j$ for $j\ge2.$
It follows that 
$1+\mathfrak{N}_p$  contains exactly $p^{2p-1}$ elements of exponent $p$, so $1+\mathfrak{N}_p\cong \Z/p^2\Z\times(\Z/p\Z)^{2p-2}$ and it is not a square of a group.
}\end{example}

As to the case $p=2$, the following example and Remark
  \ref{esempiop=2} show that the exact sequence \eqref{sep} and
  Theorem \ref{2a3b} are not sufficient to describe the 2-component of
  $A^*$.

\begin{example}
\label{nonsplit}
{\rm
Let $A=\Z[i][x,y]/(x^2-y-1, (1+i)y,y^3)$. Clearly $A$ is a type 2 ring,  $\mathfrak{N}=\mathfrak{N}_2=(\bar y)$ and as a $\Z[i]$-module 
$$
\mathfrak{N}\cong \langle\bar y\rangle\oplus \langle\bar y^2\rangle\cong (\Z[i]/(1+i))^2.
$$
Therefore, $|1+\mathfrak{N}|=|\mathfrak{N}|=2^2$.
By direct computation, one sees that $1+\bar y$ has order 4, hence the group  $1+\mathfrak{N}$ is isomorphic to $\Z/4\Z$.

Consider now the quotient ring 
$B=A/\mathfrak{N}$; we have $B=\Z[i,x]/(x^2-1)$ and 
$$B^*=B^*_2=\{\pm1,\pm i,\pm \bar x,\pm i\bar x\}\cong \Z/4\Z\times\Z/2\Z,$$
in particular $\epsilon(B)=1.$
The exact sequence \eqref{success}, that in this case coincides with \eqref{sep}, becomes
$$1\to\Z/4\Z\to A^*\to \Z/4\Z\times\Z/2\Z\to 1$$
and this gives $|A^*|=2^5.$ Finally, since $\bar x$ is a unit of order 8, one gets
$A^*\cong\langle\bar x\rangle\times\langle i\rangle\cong\Z/8\Z\times\Z/4\Z.$ 
}
\end{example}

.

\begin{remark}
\label{esempiop=2}
{\rm 
Example \ref{nonsplit} shows that 
\begin{itemize}
\item the list of 2-groups given in Proposition \ref{H2} is not exhaustive, since it does not contain the group $A^*\cong\Z/8\Z\times\Z/4\Z$;
\item it may happen that $\epsilon(A)=2$ while $\epsilon(A/\mathfrak{N})=1$;
\item the exact sequence \eqref{sep} does not split in general for $p=2$.
%; actually, it does not split in general also for $p=3$, as shown in the following example.
\end{itemize}
}
\end{remark}
%%A completely similar example of a non-splitting sequence \eqref{sep} can be made for $p=3$.
%%, which is a special case of $p\equiv 3 \pmod 4$
%\begin{example}
%\label{esempiop=3}
%{\rm
%  \label{nonsplitp=3} Take $A=\Z[i][x,y]/(x^3-y-1,3y,y^2)$. Then
%  $1+\mathfrak{N}\cong \Z/3\Z \times \Z/3\Z$ and $B=A/\mathfrak{N}=
%  \Z[i][x]/(x^3-1)$ have no elements of order 9, while $\bar x$ is a
%  unit of order 9. }
%\end{example}

\subsection{Concluding remarks}
\label{sec:Concluding remarks}

To conclude, while there is  a complete classification of the possible finite abelian groups $A^*$ when $\epsilon(A)=1$, in the case when $\epsilon(A)=2$ we have found some necessary and some sufficient conditions, which are indeed very strict.

We summarize the results for the case when  $A$ is a ring of type 2 with $\epsilon(A)=2$ as follows:
%\begin{corollary}
%\label{riassunto}
%{\sl
%Let $A$ be a ring of type 2 with $\epsilon(A)=2$. Then
\begin{itemize}
\item If $p\equiv1\pmod4$ then $(A^*)_p=1+\mathfrak{N}_p$ and can be any abelian $p$-group.
%\item If  $p\equiv3\pmod4$ and $p\ne3$ then $(A^*)_p=1+\mathfrak{N}_p$ but it can not be any $p$-group in in view of the condition on the filtration given in Proposition \ref{quadrato}.  In particular, the cardinality of  $(A^*)_p$ must be a square and all square of a $p$-group are realizable (see Proposition \ref{H2}).
%\item If $p=3$ then we have an exact sequence
%$$1\to 1+\mathfrak{N}_3\to (A^*)_3\to (\Z/3\Z)^c\to 1$$
%where $1+\mathfrak{N}_3$ is a square and $c=0$ if $\epsilon(A/\mathfrak{N})=2.$
%
\item If  $p\equiv3\pmod4$  then $(A^*)_p=1+\mathfrak{N}_p$ but it can not be any $p$-group in in view of the condition on the filtration given in Proposition \ref{quadrato}.  In particular, the cardinality of  $1+\mathfrak{N}_p$ must be a square and all squares of a $p$-group are realizable (see Proposition \ref{H2}). 
%Moreover,
%\begin{itemize}
%\item if $p>3$ then  $(A^*)_p=1+\mathfrak{N}_p$;
%\item if $p=3$ then we have an exact sequence
%$$1\to 1+\mathfrak{N}_3\to (A^*)_3\to (\Z/3\Z)^c\to 1$$
% and $c=0$ if $\epsilon(A/\mathfrak{N})=2.$
%\end{itemize}
\item If $p=2$  then we have an exact sequence
$$1\to 1+\mathfrak{N}_2\to (A^*)_2\to (\Z/2\Z)^a\times(\Z/4\Z)^b\to 1$$
where $a+b\ge1$. Moreover, Proposition \ref{H2} gives a list (not exhaustive) of realizable 2-Sylow subgroups of $A^*$. 
%and \new{ $1+\mathfrak{N}_2$ can be any abelian 2-group without direct summands of order 2.}
\end{itemize}
%}
%\end{corollary}

\smallskip

Since our description is not complete, we have tried to guess what could be a complete classification. For  the case $p=2$, what we have just said shows that the situation is far from clear, so we do not know what to expect.
The case $p\equiv3\pmod4$ appears simpler; since we have not found any example of a $p$-group with a filtration as in Proposition \ref{quadrato} which is not realizable as $(A^*)_p$, one might conjecture that the condition on the filtration could be sufficient for realizability.

\smallskip

By virtue of the decomposition $A=A_1\oplus A_2$ given in Proposition \ref{ps}, the group of units of a general ring $A$ is the product of the groups of units of a finite ring and of a type 2 ring. Therefore, the results of this section together with those of \cite{DDcharp}  allow to give a fairly precise  description of $A^*$ for a general ring $A$:
%Our results on finite characteristic rings in \cite{DDcharp} and the results of this paper 
in fact, we are able to 
produce families of new realizable abelian groups which can not be written as a product of cyclic realizable factors (already known  as a consequence of \cite{PearsonSchneider70}); on the other hand, our constraints allow to show that many families of abelian groups can not be realized.
%\smallskip
%To compare our paper to the previous work relative to finite abelian groups, we note that
%the classification of the realizable groups $A^*$ had only been made
%for cyclic groups (see \cite{Gilmer63} and \cite{PearsonSchneider70});
%of course this yields immediately that all products of these groups are realizable. 
%
%We improve upon this in two aspects: on the one hand,  our results produce families of new realizable abelian groups which can not be written as a product of cyclic realizable factors; on the other hand, our constraints allow to show that many families of abelian groups can not be realized.
%Families of examples of each type can be easily constructed from the examples already given in the paper: 
We give below one instance of each type.

\begin{example}
{\rm 
The group $\Z/4\Z\times (\Z/(11\Z))^2$ is realizable (see Proposition \ref{H2}), although $\Z/4\Z\times \Z/11\Z$ and $\Z/11\Z$ are not. 
}
\end{example}

\begin{example} 
{\rm The group $G=\Z/4\Z\times \Z/4\Z\times \Z/11\Z$
  is not realizable.  In fact, assume by contradiction that $G$ is the
  group of units of a ring $A=A_1\oplus A_2$ as in the decomposition
  by Pearson and Schneider. Then its cyclic factor $\Z/11\Z$ must come
  either from $A_1$ or from $A_2$.  It cannot come from $A_1$:
  in fact, if $11\mid |A_1^*|$ then either $10\mid |A_1^*|$ or $A_1^*$ is cyclic of order $p^\lambda -1$ 
  for suitable $p$ and $\lambda$ (see \cite[Theorem 3.1]{DDcharp}) and both these cases are impossible since $A_1^*$ is a subgroup of $G$.
   It cannot come from $A_2$ either, because in this case we
  would have $\epsilon(A_2)=2$, and, by Proposition \ref{quadrato}  the
  $11$-Sylow subgroup of $A_2^*$ must have a cardinality which is a
  square.  So $G$ is not realizable.
}
\end{example}

\smallskip

Proposition \ref{ps}, \cite[Cor. 3.2]{DDcharp} and Theorem \ref{2a3b} immediately give the following corollary  which classifies the groups of units of  reduced rings (see also Remark \ref{5.7} (ii)).
\begin{corollary}
\label{reduced}
{\sl
The finite abelian groups which are groups of units of a reduced ring are the finite products of multiplicative groups of finite fields  and possibly a group of the form  
$(\Z/2\Z)^a\times (\Z/4\Z)^b\times (\Z/3\Z)^c$
where $a,b,c\in\N$, $a+b\ge1$ and $a\ge1$ if $c\ge1.$ 
}
\end{corollary}

\medskip

Finally, we note that  combining  the previous results with \cite[Cor. 4.2]{DDcharp} we easily reobtain, as a particular case of our results, the complete classification of the finite cyclic groups which occur as group of units of a ring (firstly given in \cite{PearsonSchneider70}). 

\begin{corollary}
\label{ciclici}
{\sl
A finite cyclic group is the group of units of ring  
if and only if its order is  the product of a set of pairwise coprime integers of the following list:
\begin{enumerate}[label={\alph*)}]
\item $p^\lambda-1$ where $p$ is a prime and $\lambda\ge1$;
\item $(p-1)p^k$  where $p>2$ is a prime and $k\ge1$;
\item $2d$ where $d>0$ is odd; 
\item  $4d$ where $d$ is an odd integer and each of its prime factors is congruent to 1 mod 4.
 \end{enumerate}
 }
 \end{corollary}
\begin{proof}
By \cite[Cor. 4.2]{DDcharp}, the cyclic groups of units which are realizable with finite characteristic ring are those whose order is  the product of pairwise coprime integers as in $(a)$, $(b)$ and, possibly, 2 or 4.
Consider now rings of characteristic 0: by Theorem \ref{2n} all cardinalites listed in $(c)$ are possible; by Proposition \ref{H2} and Corollary \ref{nocyclic} the cardinality $4d$ is possible if and only if $d$ is an odd integer and each of its prime factors is congruent to 1 mod 4.
\end{proof}

\medskip

	\appendix{}
	
\section{Densities}
\label{densities}

In this section we study the distribution of the possible values of
the cardinality of $A^*$ where $A^*$ is a finite abelian group. The
first simple remark regards the prime numbers $p$ that occur as the
cardinality of $A^*$: by Corollary \ref{corcor}, we know that the only possibilities are $p=2$ or $p$ a Mersenne prime.
%By Theorem \ref{2n}, in characteristic zero the only possibility is
%$p=2$, while, by \cite[Corollary 3.4]{DDcharp}, in finite characteristic we can have odd primes
%$p$ if and only if they are of type $p=2^k-1$ for some $k>1$. In conclusion, we have the following
%
%
%
%\begin{prop}
%  \label{mersenne} {\sl The possible prime
%    values of $|A^*|$ are $p=2$ and the Mersenne primes.}
%\end{prop}

\smallskip

The next observation is that odd numbers can be very seldom the cardinality of $A^*$ for some ring $A$. 
For a set $X\subset \N$ and for $n\in \N$ we denote, as usual, $X(n):=X\cap\{1,\dots,n\}$ and 
the density of $X$ as the limit
$$ 
\lim_{n\to\infty} \frac {|X(n)|}{n} 
$$ 
whenever this limit exists. We have: 

\begin{prop}
\label{ddispari}
{\sl The set of all possible odd values of $|A^*|$ has density zero.}
\end{prop}
\begin{proof}
  Let $X=\{\, |A^*|\, :\, |A^*|\ \hbox{is odd}\,\}$. By
  \cite[Corollary 3.4]{DDcharp} and Theorem \ref{2n},
$X=\{\prod_{i=1}^m (2^{k_i}-1)\mid m\in\N,\ k_i\ge2\}$. For any
positive number $M$ we estimate the cardinality of the set $X(2^M)=\{
n\in X \mid n\le 2^M\}$. Clearly,
$$2^{\sum_{i=1}^m(k_i-1)}=\prod_{i=1}^m2^{k_i-1}\le\prod_{i=1}^m(2^{k_i}-1),$$
and if $n=\prod_{i=1}^m(2^{k_i}-1)\in X(2^M)$, then $\sum_{i=1}^m(k_i-1)\le M$.
It follows that $|X(2^M)|$ is at most $\sum_{h=1}^M p(h)$, where $p(h)$ denotes the number of partitions of the positive integer $h$. By \cite[Thm 6.3]{Andrews84}
$$p(h)\sim \frac1{4h\sqrt3}\exp\left(\pi\left(\small{\frac23}\right)^\frac12 h^\frac12\right),$$
so 
$$\frac1{2^M}\sum_{h=1}^M p(h)\ll \frac1{2^M}\log M\exp\left(\pi\left(\frac {\small 2}{\small 3}\right)^\frac12 M^\frac12\right)$$
  and the last expression goes to zero as $M\to\infty.$
\end{proof}

Theorem \ref{2n} says that any even number can be the cardinality of $A^*$ for some ring $A$. Together with Proposition \ref{ddispari} this gives

\begin{prop}
\label{densitatotale}
{\sl The density of the possible cardinalities of $A^*$ is equal to $\frac 12$.} 
\end{prop}

As we have seen before, the possible groups $A^*$ are much less if we restrict to the case where the ring $A$ has no non-zero nilpotents. In fact, we have the following

\begin{prop}
\label{density-nonilp}
{\sl The set of possible cardinalities if $A^*$, where $A$ ia a reduced ring, has density zero. }
\end{prop}

\begin{proof}
We need the following 

\begin{lemma} 
\label{densita2^k}
{\sl For each $h\ge 0$, the set $X_h$ of possible cardinalities of $A^*$ such that $2^h||{\rm card}(A^*)$, where $A$ is a reduced ring, has density zero.}
\end{lemma}

\begin{proof}
  For $h=0$, the statement has already been proved in Proposition
  \ref{ddispari}. For $h>0$, observe that, by Corollary \ref{reduced},  the elements of $m\in X_h$ are
  obtained as products of type 
\begin{equation}
\label{prodotti}
m = 2^{a+2b}\prod_{i=1}^s (2^{k_i}-1) \prod_{j=1}^t (p_j^{\lambda_j}-1)\, , 
\end{equation} 
where $a+2b\le h$, $p_1,\dots,p_s$ are odd primes and $2^{h-a-2b} ||
\prod_{j=1}^t (p_j^{\lambda_j}-1)$ (the possible  factor
$3^c$ given by Theorem \ref{2a3b} can be absorbed in the product $\prod_{i=1}^s (2^{k_i}-1)$, since
$2^2-1=3$).

\smallskip

We write a number $m\in X_h(n)$ as $m=2^{a+2b}m_1m_2$, where $m_1=\prod_{i=1}^s (2^{k_i}-1)$ and $m_2=\prod_{j=1}^t (p_j^{\lambda_j}-1)$ and we partition $X_h(n)$ according to the value of $m_1\le n$.
We have $|X_h(n)|\le \sum_{m_1\le n} |\{2^{a+2b}m_2\le n/m_1\}|$.  To estimate the cardinality of each summand we first note that  $m_2\le \frac n{m_1}$ and, since all
numbers $p^{\lambda_j}-1$ are even, we have $t\le h$. Moreover,  
we have trivially
$p^{\lambda}\le\frac 32(p^\lambda-1)$, so $\prod_{j=1}^t p^{\lambda_j}
\le (\frac 32)^t \frac n{m_1} \le (\frac 32)^h \frac n{m_1}$. 

It follows that, for any $m_1\le n$ and any given pair $(a,b)$ with
$a+2b\le h$, the number of possible $m_2$ for which $m\in X_h(n)$
does not exceed the cardinality of the set $\{\prod_{j=1}^t p^{\lambda_j}\le (\frac 32)^h \frac
n{m_1}\}$.  
For any positive real number $x$ and any positive integer $t$, let
$$
\omega_t(x) = \left|\{p_1^{\lambda_1}\dots p_t^{\lambda_t} \le x\mid \lambda_j\ge0\}\right|\, .
$$
Clearly $\omega_t(x)\le \omega_h(x)$ for $t\le h.$ 
By \cite[Lemma B]{Ramanujan17}, we have 
$$
\omega_h(x)=O\left(\frac {x (\log\log x)^{h-1}}{\log x}\right).
$$
Putting $x= (\frac 32)^h \frac n{m_1}$ and summing over the finitely many possibilities for the pairs $(a,b)$, we obtain that 

$$
|\{2^{a+2b}m_2\le n/m_1\}| \le c(h)\,  \omega_h\left(\left(\frac 32\right)^h \frac n{m_1}\right)
$$
for some positive constant $c(h)$ depending only on $h$. 
Next, we split the numbers in $X_h(n)$ into two subsets, according to the value of $m_1$. We have 
\begin{equation}
\label{twosums}
|X_h(n)| \le c(h)\sum_{m_1\le \sqrt{(\frac 32)^hn}} 
\omega_h \left(\frac{(\frac 32)^hn}{m_1}\right) + c(h)\sum_{m_1>\sqrt{(\frac 32)^hn}} 
\omega_h\left(\frac{(\frac 32)^hn}{m_1}\right)\, .
\end{equation}

% \begin{equation}
% \label{twosums}
% \begin{split}
% |X_h(n)| \le \sum_{\genfrac{}{}{0pt}{}{k_1,\dots,k_s}{(2^{k_1}-1)\cdots(2^{k_s}-1)\le \sqrt{(\frac 32)^hn}}} 
% \omega_h \left(\frac{(\frac 32)^hn}{(2^{k_1}-1)\cdots(2^{k_s}-1)}\right) \\ 
%  + \sum_{\genfrac{}{}{0pt}{}{k_1,\dots,k_s}{(2^{k_1}-1)\cdots(2^{k_s}-1)>\sqrt{(\frac 32)^hn}}} \omega_h\left(\frac{(\frac 32)^hn}{(2^{k_1}-1)\cdots(2^{k_s}-1)}\right)\, .
% \end{split}
% \end{equation}
For a generic term of the first sum in \eqref{twosums}  we have
$ \log \left(\frac{(\frac 32)^hn}{m_1}\right)\gg \log n$, so, for $n$ sufficiently large,  the contribution of each of these terms is less
than or equal to $c_1(h)\,\frac1{m_1}\frac {n\, (\log\log n)^{h-1}}{\log n}$ for some
  positive constant $c_1(h)$.  Now, 
$$
\sum_{m_1\le n} \frac 1{m_1} \le \sum_{k_1,\dots,k_s} \frac 1{(2^{k_1}-1)\cdots(2^{k_s}-1)} :=S
$$ 
and the series $S$ converges, since 
$$ S\le \prod_{k=1}^\infty \frac {2^k-1}{2^k-2} = \prod_{k=1}^\infty \left(1+\frac 1{2^k-2}\right) $$ 
and the series of logarithms $\sum_{k=1}^\infty \log(1+\frac1{2^k-2})\sim \sum_{k=1}^\infty \frac 1{2^k-2}$ converges. It follows that there is a positive constant $c_2$ such that the first sum is less than or equal to 
$c_2 \, \frac {n(\log\log n)^{h-1}}{\log n}$ and therefore is $o(n)$.

As to the second sum, each term in it is clearly $\ll n^{\frac 12}$, so we are reduced to estimate the number of terms is this sum. 

We now write the number $m_1$ as 
$(2^2-1)^{\mu_1}(2^3-1)^{\mu_2}\cdots(2^{r+1}-1)^{\mu_{r}}$; since $m_1\le n$, then {\it a fortiori},
$2^{\mu_1}2^{2\mu_2}\cdots 2^{r\mu_r}\le n$, so the non-negative integers
$\mu_i$ are such that $\mu_1+2\mu_2+\dots+r\mu_r\le\frac{ \log n}{\log 2}$. It
is well-known that the number of integer points $(\mu_1,\dots,\mu_r)$
satisfying the preceding inequality is approximately equal to the
volume of the region
$$ \{(u_1,\dots,u_r)\in \R^r\mid u_i\ge 0, \ u_1+2u_2+\dots +ru_r\le \frac{ \log n}{\log 2}\}.$$ 
More precisely, by \cite[Chapter 6, Thm 2]{Lang94ANT}, we get that the number of
$r$-tuples $(\mu_1,\dots,\mu_r)$ such that $2^{\mu_1}2^{2\mu_2}\cdots 2^{r\mu_r}\le n$ is less than 
$$
\frac 1{(r!)^2}\frac{(\log n)^r}{(\log 2)^r}+O\left(r\, \frac{(\log n)^{r-1}}{(r-1)!^2})\right)\, .
$$
By the trivial inequality $(2r)!\ge (r!)^2$ we get 
\begin{eqnarray*}
\frac {(\log n)^r}{(r!)^2}&\le& \frac{(\log n)^r}{(2r)!} \\
&=& \frac {((\log n)^{\frac 12})^{2r}}{(2r)!} \\
&\le&  e^{(\log n)^{\frac 12}} \\
&\ll& (\log n) e^{(\log n)^{\frac 12}}  = o(n^\epsilon)
\end{eqnarray*}
for any $\epsilon >0$. Taking into account that $r\ll \log n$, the error term is estimated similarly, so the second sum is also $o(n)$. 

Putting together the last two estimates we get the lemma. 

\end{proof}

Denote by $X'_h$ the set of all possible cardinalities of $A^*$
divisible by $2^h$, where $A$ is a reduced ring. For any $h\ge 0$, we
can obviously partition the set $X$ of all possible cardinalities of $A^*$ as
$$ X = X_0\cup X_1\ \cup \dots \cup X_{h-1}\cup X'_h. $$

By the lemma, the densities $\delta(X_i)$ of $X_i$ are zero for
$i=0,\dots,h=1$, so the density of $X$ equals the density of $X'_h$. 
But trivially $\delta(X'_h)\le \frac 1{2^{h+1}}$, so, taking the limit for
$h\to\infty$, we get $\delta(X)=0$.
\end{proof}

% \begin{prop}
% \label{proddir}
% {\sl Let $A$ be a ring of characteristic $n>0$, and let $n=p_1^{a_1}\cdots p_r^{a_r}$ be the factorization of $n$. Then 
% $$ A \cong A_1\times \dots, \times A_r,$$
% where ${\rm char} (A_i) = p_i^{r_i}$ for $i=1,\dots,r$. 
% }
% \end{prop}

% \begin{corollary}
% \label{descrizioneAi}
% {\sl In the preceding notation, $A_i\cong \prod\limits_{j\ne i} p_j^{a_j}A$.} 
% \end{corollary}

\bibliographystyle{amsalpha}
\bibliography{biblio}

% \bibliography{mybib}{}
% \bibliographystyle{plain}

% \begin{thebibliography}{99}
% \bibitem{Andrews}  \textsc {G. E. Andrews},
%                 \textit {The Theory of Partitions},
%                \textup {Cambridge University Press 1984}.
               
%  \bibitem{ditor} \textsc{S. Z. Ditor}, On the group of units of a ring,                
%  \textit{Amer. Math. Monthly} {\bf 78} (1971), 522-523.

%  \end{thebibliography}
\end{document}